\documentclass[a4paper, leqno, 12pt]{article}
\usepackage{amsmath,amsthm}
\usepackage{amssymb,latexsym}
\usepackage{enumerate}
\usepackage{color}
\usepackage[T1]{fontenc}
\usepackage{array}
\usepackage[all]{xypic}

\usepackage[left=3cm,right=3cm,top=3.0cm,bottom=3.0cm]{geometry}

\newtheorem{thm}{Theorem}[section]
\newtheorem{cor}[thm]{Corollary}

\newtheorem{prop}[thm]{Proposition}

\theoremstyle{definition}

\numberwithin{equation}{section}

 1

\DeclareMathOperator{\N}{\mathbb{N}}

\DeclareMathOperator{\ra}{\rightarrow}

\DeclareMathOperator{\lra}{\leftrightarrow}
 
\DeclareMathOperator{\ran}{\rangle}\DeclareMathOperator{\lan}{\langle}
\DeclareMathOperator{\para}{\parallel}\DeclareMathOperator{\npara}{\nparallel}

\def\NN{{\mathbb{N}}}
\def\ZZ{{\mathbb{Z}}}

\def\FF{{\mathbb{F}}}

\def\CB{{\cal B}}

\def\CF{{\cal F}}

\def\CL{{\cal L}}

\def\vsp{\vspace*{1.5ex}}

\def\epv {{$\mbox{}$\hfill ${\Box}$\vspace*{1.5ex} }}
\def\rk{\textnormal{rk}}

\def\ra{\rightarrow}

\def\und#1{\underline{#1}}

\def\ACF{\textnormal{ACF}} \def\TF{\textnormal{TF}}

\def\CIS{\textnormal{CIS}}

\begin{document}


\baselineskip=17pt


\title{A constructive proof of Tarski's theorem on quantifier elimination in the theory of ACF}
\author{Grzegorz Pastuszak (Toru\'n)}

\date{}

\maketitle

\renewcommand{\thefootnote}{}

\renewcommand{\thefootnote}{\arabic{footnote}}
\setcounter{footnote}{0}


\begin{abstract} Assume that $\ACF$ denotes the theory of algebraically closed fields. The renowned theorem of A. Tarski states that $\ACF$ admits quantifier elimination. In this paper we give a constructive proof of Tarski's theorem on quantifier elimination in $\ACF$. This means that for a given formula $\varphi$ of the language of fields we construct a quantifier-free formula $\varphi'$ such that $\ACF\vdash\varphi\lra\varphi'$. We devote the last section of the paper to show some applications of this constructive version in mathematics and physics.

\end{abstract}

\section{Introduction and notation}

Throughout the paper, $\CF=(0,1,+,-,\cdot)$ denotes the language of fields, $\TF$ the theory of fields and $\ACF$ the theory of algebraically closed fields. We use the notation and terminology of \cite{M} for the basic concepts of model theory. In most cases, our notation is the standard one. Nevertheless, we recall some notions at the end of this section for convenience of the reader.

In 1948 A. Tarski proved (in an unpublished paper, see \cite{Rob} for the details) that $\ACF$ admits quantifier elimination. This is one of the most fundametal facts in model theory. Therefore, there is a number of proofs of Tarski's theorem in the literature. Standard ones are existential, that is, they do not provide the form of the quantifier-free formula equivalent with the given one. This paper aims to provide that form. More precisely, for a given $\CF$-formula $\varphi$ we construct a quantifier-free formula $\varphi'$ such that $\ACF\vdash\varphi\lra\varphi'$. Our construction is based on results of \cite{P}. In that paper we set a bound on the length of ascending chains of ideals in multivariate polynomial rings. These ideals are generated by polynomials of degrees less or equal to fixed natural numbers. In a sense, we rediscovered in \cite{P} some of the main results of \cite{Mo} and \cite{AsPo} (see also \cite{Se1}) in order to prove Tarski's theorem in a constructive way.

We emphasize that the results of \cite{Mo} and \cite{AsPo}, together with argumentation similar to that of Section 3, enable to give an alternative constructive proof of Tarski's theorem. Moreover, there are many results on effective quantifier elimination in $\ACF$, see for example \cite{PS}. Therefore we are not pioneer in these considerations.

The paper is organized as follows. In Section 2 we recall some results from \cite{P}. We aim to present Corollary 2.3 (Corollary 4.5 in \cite{P}) which is the main tool in our proof. The constructive proof of Tarski's theorem is presented in Section 3 which is the core of the paper. The main result is Theorem 3.2. As a corollary of Theorem 3.2 we get a computable condition for the existence of a common root of multivariate polynomials, see Corollary 3.3. By \textit{computable condition} (or \textit{computable criterion}) we mean a procedure employing only finite number of arithmetic operations. The last section of the paper is devoted to show examples of application of Theorem 3.2 in mathematics and physics. These applications are connected with some problems of quantum information theory which we consider in \cite{JP} and \cite{PJ} (see also \cite{JKP} and \cite{PKJ} for similar topics).

The results of Section 3 are part of the author's master's thesis, supervised by Stanis{\l}aw Kasjan in 2007. The author is grateful to the supervisor for all discussions and support during the work on the thesis.

We introduce some notation and terminology. Assume that $\CL$ is a language and $\varphi_{1},...,\varphi_{n}$ are $\CL$-formulas. Then $\bigwedge_{i=1}^{n}\varphi_{i}$ and $\bigvee_{i=1}^{n}\varphi_{i}$ denote the formulas $\varphi_{1}\wedge...\wedge\varphi_{n}$ and $\varphi_{1}\vee...\vee\varphi_{n}$, respectively. If $\und{x}=(x_{1},...,x_{m})$ is a sequence of variables and $Q$ is a quantifier, then $Q_{\und{x}}$ is the abbreviation of $Q_{x_{1}}...Q_{x_{m}}$. Generally, if $A=\{a_{1},...,a_{s}\}$ is a set of variables, then $Q_{A}$ is the abbreviation of $Q_{b_{1}}...Q_{b_{s}}$ where $b_{1},...,b_{s}$ is any permutation of $a_{1},...,a_{s}$. This is consistent since, for any $\CL$-formula $\varphi$, the formulas $Q_{b_{1}}...Q_{b_{s}}\varphi$ and $Q_{a_{1}}...Q_{a_{s}}\varphi$ are equivalent. 

If $\varphi$ is an $\CL$-formula and $a_{1},...,a_{n}$ are all free variables of $\varphi$, then sometimes we write $\varphi(\und{a})$ instead of $\varphi$ where $\und{a}=(a_{1},...,a_{n})$. Recall that if $\varphi(\und{a})$ is an atomic $\CF$-formula ($\CF$ denotes the language of fields), then $\varphi(\und{a})$ has the form $F=0$ or $F\neq 0$ where $F$ is a mutivariate polynomial in $\ZZ[a_{1},...,a_{n}]$. 

If $\und{a}=(a_{1},...,a_{n})$, then $\ZZ[\und{a}]$ denotes the ring $\ZZ[a_{1},...,a_{n}]$. If $\und{x}=(x_{1},...,x_{m})$, then $\ZZ[\und{a}][\und{x}]$ denotes the ring of polynomials in $m$ variables $x_{1},...,x_{m}$ over the ring $\ZZ[\und{a}]$. A polynomial $F$ in $\ZZ[\und{a}][\und{x}]$ has the form $\sum_{\alpha\in\NN^{m}}f_{\alpha}\cdot\und{x}^{\alpha}$ where $f_{\alpha}\in\ZZ[\und{a}]$ for any $\alpha\in\NN^{m}$ and $f_{\alpha}=0$ for almost all $\alpha\in\NN^{m}$. Here, $\und{x}^{\alpha}$ denotes $x_{1}^{\alpha_{1}}...x_{m}^{\alpha_{m}}$ where $\alpha=(\alpha_{1},...,\alpha_{m})\in\NN^{m}$. The degree of $F$ with respect to $x_{1},...,x_{m}$ is denoted by $\deg(F)$. Generally, if $C$ is a set of variables, then $\ZZ[C][\und{x}]$ is the ring of polynomials in $m$ variables $x_{1},...,x_{m}$ over the ring $\ZZ[C]$ of polynomials in variables from $C$.

We denote by $\NN$ the set of all natural numbers and by $\NN_{1}$ the set $\NN\setminus\{0\}$. Assume that $m\in\NN_{1}$. We view the set $\NN^{m}$ as a monoid with respect to the pointwise addition, denoted by $+$. We denote by $\und{0}$ the neutral element $(0,...,0)\in\NN^{m}$ of $+$. If $\alpha,\beta\in\NN^{m}$ and $\alpha+\gamma=\beta$ for some $\gamma\in\NN^{m}$, then we write $\alpha\para\beta$. Note that $\para$ defines an order on $\NN^{m}$ and $\NN^{m}$ is an ordered monoid with respect to $+$ and $\para$. Sometimes we treat the elements of the set $\NN^{m}$ as sequences of natural numbers. If $\alpha\in\NN^{m}$ and $\alpha=(a_{1},...,a_{m})$, then we set $|\alpha|=a_{1}+...+a_{m}$. 

\section{Ascending chains of ideals in the polynomial ring}

In this section we recall the results of \cite{P} which are the main tool in constructive proof of Tarski's theorem. The first goal is to recall the construction of a function with the \textit{bounding property}. We use this function in Theorem 2.2 to set a bound on the length of ascending chains of ideals in $K[x_{1},...,x_{m}]$ ($K$ is a field) which are generated by polynomials of degrees less or equal to fixed natural numbers. Then we present Corollary 2.3 which we directly apply in the proof of Tarski's theorem. This section does not contain any proof. We refer to \cite{P} for all the proofs and other details.


We denote by $\FF$ the set of all non-decreasing functions $\N_{1}\ra\N_{1}$. We write $f\leq f'$ if and only if $f,f'\in\FF$ and $f(n)\leq f'(n)$ for any $n\in\NN_{1}$. If $f\in\FF$ and $s\in\NN$, then ${}^{s}f:\NN_{1}\ra\NN_{1}$ is a function such that ${}^{s}f(n)=f(s+n)$ for any $n\in\NN_{1}$. Observe that ${}^{s}f\in\FF$. A sequence $\alpha_{1},...,\alpha_{t}\in\NN^{m}$ is an \textit{antichain} if and only if $\alpha_{i}\npara\alpha_{j}$ for any $i<j$. Assume that $f\in\FF$. We say that an antichain $\alpha_{1},...,\alpha_{t}\in\NN^{m}$ is \textit{$f$-bounded} if and only if $|\alpha_{i}|\leq f(i)$ for any $i=1,...,t$. Assume that $m\geq 1$ is a natural number. We say that a function $\CB_{m}:\FF\ra\NN$ has the \textit{bounding property for $m$} if and only if the following conditions are satisfied:
\begin{enumerate}[\rm(1)]
	\item $t\leq\CB_{m}(f)$ for any $f\in\FF$ and $f$-bounded antichain $\alpha_{1},...,\alpha_{t}\in\NN^{m}$ of length $t$,
	\item $\CB_{m}(f)\leq\CB_{m}(f')$ for any $f,f'\in\FF$ such that $f\leq f'$.
\end{enumerate}


We say that a function $\CB:\N_{1}\times\FF\ra\NN$ has the \textit{bounding property} if and only if, for any $m\in\N_{1}$, the function $\CB_{m}:\FF\ra\NN$ defined by $\CB_{m}(f)=\CB(m,f)$, for any $f\in\FF$, has the bounding property for $m$.  

The existence of a function with the bounding property is a consequence of the Compactness Theorem of first order logic, see \cite{F} and \cite[Proposition 3.25]{AsPo} for more details. However, this approach does not provide the explicit form of such a function.

We recall from \cite{P} the construction of a function with the bounding property. Equivalently, we give a sequence $(\CB_{m})_{m\in\NN_{1}}$ of functions such that $\CB_{m}:\FF\ra\NN$ has the bounding property for $m$. The construction is inductive with respect to the number $m$. It is given in two main steps, but the second step is divided in three parts.

Step 1. Assume that $m=1$. We define $\CB_{1}:\FF\ra\NN$ to be a function such that $\CB_{1}(f)=f(1)+1$ for any $f\in\FF$.

Step 2. Assume that $m\geq 2$ and the function $\CB_{m-1}:\FF\ra\NN$ is defined. In order to define $\CB_{m}:\FF\ra\NN$, we construct some sequence of functions $(\CB_{m}^{k})_{k=0}^{m}$, $\CB_{m}^{k}:\FF\times\NN^{k}\ra\NN$. This is done by the backward induction with respect to the number $k$. We give the construction in three steps.

Step 2.1. Assume that $k=m$. We define $\CB_{m}^{m}:\FF\times\NN^{m}\ra\NN$ to be a function such that $\CB_{m}^{m}(f,b_{1},...,b_{m})=(b_{1}+1)\cdot...\cdot(b_{m}+1)$ for any $f\in\FF$ and $(b_{1},...,b_{m})\in\NN^{m}$.

Step 2.2. Assume that $k\in\{0,...,m-1\}$ and the function $\CB_{m}^{k+1}:\FF\times\NN^{k+1}\ra\NN$ is defined. Suppose $f\in\FF$, $\beta\in\NN^{k}$ and let $g:\NN_{1}\ra\NN_{1}$ be a function such that $g(1)=1$ and $$g(n+1)=1+g(n)+\CB^{k+1}_{m}({}^{g(n)}f,\beta,f(g(n)))$$ for any $n\geq 1$. We have $g\in\FF$ and hence there is a function $\CF_{m}^{k}:\FF\times\NN^{k}\ra\FF$ such that $(f,\beta)\mapsto g$. We set $\CB_{m}^{k}(f,\beta)=g(\CB_{m-1}(f\circ g)+1)$ for any $f\in\FF$, $\beta\in\NN^{k}$ and $g=\CF_{m}^{k}(f,\beta)$.

Step 2.3. We identify $\CB_{m}$ with $\CB_{m}^{0}$.

\vsp
The above procedure defines a sequence of functions $(\CB_{m})_{m\in\NN_{1}}$, $\CB_{m}:\FF\ra\NN$. Let $\CB:\NN_{1}\times\FF\ra\NN$ be a function such that $\CB(m,f)=\CB_{m}(f)$ for any $m\in\NN_{1}$ and $f\in\FF$. In Section 3 of \cite{P} we prove the following theorem.

\begin{thm} The function $\CB_{m}:\FF\ra\NN$ has the bounding property for $m$, for any $m\in\NN_{1}$. Consequently, the function $\CB:\NN_{1}\times\FF\ra\NN$ has the bounding property.
\end{thm}

{\bf Proof.} See Proposition 3.1, Proposition 3.2 and Corollary 3.4 from \cite{P}. \epv

Assume that $K$ is a field, $m\geq 1$ is a natural number and $f:\N_{1}\ra\N_{1}$ is an arbitrary function. An ascending chain $I_{1}\subsetneq ...\subsetneq I_{t}$ of ideals in $K[x_{1},,...,x_{m}]$ is \textit{$f$-bounded} if and only if $I_{j}$ is generated by polynomials of degrees less or equal to $f(j)$, for any $j=1,...,t$. 

Theorem 2.1 is used in \cite{P} to give a bound on the length of $f$-bounded ascending chains of ideals in $K[x_{1},,...,x_{m}]$ depending on $m$ and $f$. We recall the appropriate theorem below.

\begin{thm} Assume that $m\geq 1$ and $f:\NN_{1}\ra\NN_{1}$ is a function. Suppose that $I_{1}\subsetneq ...\subsetneq I_{t}$ is an $f$-bounded ascending chain of ideals in $K[x_{1},,...,x_{m}]$ of length $t$. Let $g:\NN_{1}\ra\NN_{1}$ be a non-decreasing function such that $g(n)$ is the greatest number of the set $\{f(1),f(2),...,f(n)\}$, for any $n\in\NN$. Then $t\leq\CB(m,g)$. In particular, we have $t\leq\CB(m,f)$, if $f$ is non-decreasing. 
\end{thm}

{\bf Proof.} See Theorem 4.2 from \cite{P}. \epv

Let $d\geq 1$ be a fixed natural number. By a string $3^{n}d$ we mean the function $f:\NN_{1}\ra\NN_{1}$ such that $f(n)=3^{n}d$. Set $m\geq 1$, $d\geq 1$ and define the function $\gamma_{m,d}:\NN\ra\NN$ such that $$\gamma_{m,d}(i)=(3^{\CB(m,3^{n}d)-1}-1)d+i$$ for any $i\in\NN$. Applying Theorem 2.2 and the theory of Gr\"obner bases (see \cite{AL}) we prove in \cite{P} the following result which plays a crucial role in constructive proof of Tarski's theorem.

\begin{cor} Assume that $m\geq 1$ and $d\geq 1$. Then for any $G\in K[x_{1},...,x_{m}]$ and $F_{1},...,F_{s}\in K[x_{1},...,x_{m}]$ such that $\deg(F_{i})\leq d$ for $i=1,...,s$ the following condition is satisfied: $G\in\lan F_{1},...,F_{s}\ran$ if and only if there exist $H_{1},...,H_{s}\in K[x_{1},...,x_{m}]$ such that $G=H_{1}F_{1}+...+H_{s}F_{s}$ and $\deg(H_{i})\leq\gamma_{m,d}(\deg(G))$ for $i=1,...,s$.
\end{cor}

{\bf Proof.} See Proposition 4.3, Corollary 4.4 and Corollary 4.5 from \cite{P}. \epv

\section{Tarski's theorem}

This section is devoted to the constructive proof of Tarski's theorem. We recall that it is enough to give the construction for some special formulas over the language $\CF$ of fields which we call \textit{common root formulas}.

Let $\CL$ be a language and assume that $\varphi$ is a formula over $\CL$. It is well known that $\varphi$ can be written in \textit{prenex normal form}, see for example \cite[Chapter 3]{Rot}. It follows from De Morgan's laws that $\varphi$ is equivalent with the formula $\bigvee_{i=1}^{t}\exists_{\und{x}}(\bigwedge_{j=1}^{s_{i}}\varphi_{ij})$ where $\varphi_{ij}$ are atomic formulas or negations of atomic formulas. 

An $\CL$-formula is a \textit{conjunctive prenex normal formula} if it has the form $\exists_{\und{x}}(\bigwedge_{i=1}^{s}\varphi_{i})$ where each $\varphi_{i}$ is an atomic $\CL$-formula or a negation of such. Hence a theory $T$ over $\CL$ admits quantifier elimination if an only if for any conjunctive prenex normal $\CL$-formula $\varphi$ there is a quantifier-free $\CL$-formula $\varphi'$ such that $T\vdash\varphi\lra\varphi'$. We recall below the form of conjunctive prenex normal $\CF$-formulas.

Assume that $\und{a}=(a_{1},...,a_{n})$, $\und{x}=(x_{1},...,x_{m})$ and $F_{1},...,F_{s}\in\ZZ[\und{a}][\und{x}]$. Assume that $F_{i}=\sum_{\alpha\in\NN^{m}}f_{i,\alpha}\cdot\und{x}^{\alpha}$ where $f_{i,\alpha}\in\ZZ[\und{a}]$ for any $i=1,...,s$, $\alpha\in\NN^{m}$ and $f_{\alpha}=0$ for almost all $\alpha\in\NN^{m}$. A formula of the form $\exists_{\und{x}}(F_{1}(\und{x})=0\wedge...\wedge F_{s}(\und{x})=0)$ is a \textit{common root formula}. 

\begin{prop}
Any conjunctive prenex normal $\CF$-formula is equivalent with some common root formula.
\end{prop}

{\bf Proof.} Assume that $\und{a}=(a_{1},...,a_{n})$ and $\varphi(\und{a})$ is a conjunctive prenex normal $\CF$-formula. Then $$\varphi(\und{a})=\exists_{\und{x}}(F_{1}(\und{x})=0\wedge...\wedge F_{r}(\und{x})=0\wedge G_{1}(\und{x})\neq 0\wedge...\wedge G_{t}(\und{x})\neq 0)$$ where each $F_{i},G_{j}$ is a polynomial of the form $\sum_{\alpha\in\NN^{m}}f_{\alpha}\cdot\und{x}^{\alpha}$ where $f_{\alpha}\in\ZZ[\und{a}]$ and $f_{\alpha}=0$ for almost all $\alpha\in\NN^{m}$. Since the formula $G_{1}(\und{x})\neq 0\wedge...\wedge G_{t}(\und{x})\neq 0$ is equivalent with $(G_{1}\cdot...\cdot G_{t})(\und{x})\neq 0$, the formula $\varphi(\und{a})$ is quivalent with $$\varphi'(\und{a})=\exists_{\und{x},z}(F_{1}(\und{x})=0\wedge...\wedge F_{r}(\und{x})=0\wedge zG(\und{x})-1=0)$$ where $G=G_{1}\cdot...\cdot G_{t}$. This shows the assertion. \epv


Common root formulas play a crucial role in the constructive proof of Tarski's theorem. We aim to give an equivalent quantifier-free form of common root formulas. 

Assume that $d,d'\geq 1$ are some fixed natural numbers. Let $F_{1},...,F_{s}\in\ZZ[\und{a}][\und{x}]$ be polynomials such that $\deg(F_{i})\leq d$ and $F_{i}=\sum_{|\alpha|\leq d}f_{i,\alpha}\cdot\und{x}^{\alpha}$ where $f_{i,\alpha}\in\ZZ[\und{a}]$ for any $i=1,...,s$ and $\alpha\in\NN^{m}$. Let $A^{d,d'}_{F_{1},...,F_{s}}=A$ be a matrix with rows indexed by elements of the set $X=\{\delta\in\NN^{m}|d+d'\geq|\delta|\}$, columns indexed by elements of $\{1,...,s\}\times X$ and $$A(\delta,(i,\beta))=\left\{\begin{array}{cccc}f_{i,\delta-\beta}&&\textnormal{if $\beta\para\delta$,}\\0&&\textnormal{otherwise}\end{array}\right.$$ where $\delta,\beta\in X$ and $i\in\{1,...,s\}$. Let $\widehat{A}^{d,d'}_{F_{1},...,F_{s}}=\widehat{A}$ be an augmented matrix $(A|B)$ where $B$ is a column with $\{1,...,s\}\times X$ rows such that $B=\left[\begin{matrix}0&&\hdots&&0&&1\end{matrix}\right]^{T}$. Assume that $S(A)$ and $S(\widehat{A})$ are the sets of all square submatrices of $A$ and $\widehat{A}$, respectively. Moreover, assume that $S(\widehat{A},n)$ is the subset of $S(\widehat{A})$ consisting of the matrices of order greater than $n$. We define a quantifier-free formula $$\Delta^{d,d'}_{F_{1},...,F_{s}}(\und{a})=\bigwedge_{M\in S(A)}(\det M\neq 0\ra(\bigvee_{N\in S(\widehat{A},o_{M})}\det D\neq 0))$$ where $o_{M}$ denotes the order of the matrix $M$. Assuming that $\und{a}$ is a tuple of elements of some field, the formula $\Delta^{d,d'}_{F_{1},...,F_{s}}(\und{a})$ holds if and only if the rank of the matrix $\widehat{A}$ is greater then the rank of $A$.


In the following theorem we show that common root formulas are equivalent with quantifier-free formulas of the form $\Delta^{d,d'}_{F_{1},...,F_{s}}(\und{a})$. This theorem is a constructive version of Tarski's theorem on quantifier elimination in the theory of ACF, because any $\CF$-formula can be easily written as a disjunction of common root formulas.

The aforementioned theorem is the main result of the paper. The proof is based on Corollary 2.3 and hence on the results of \cite{P} recalled in Section 2. 

\begin{thm} Assume that $\und{a}=(a_{1},...,a_{n})$, $\und{x}=(x_{1},...,x_{m})$, $F_{1},...,F_{s}\in\ZZ[\und{a}][\und{x}]$ and $\varphi(\und{a})=\exists_{\und{x}}(F_{1}(\und{x})=0\wedge...\wedge F_{s}(\und{x})=0)$. Assume that $d$ is the maximum of degrees of polynomials $F_{1},...,F_{s}$ and $d'=\gamma_{m,d}(0)$. Then $\ACF\vdash\varphi(\und{a})\lra\Delta^{d,d'}_{F_{1},...,F_{s}}(\und{a})$.
\end{thm}

{\bf Proof.} Assume that $F_{i}=\sum_{|\alpha|\leq d}f_{i,\alpha}\cdot\und{x}^{\alpha}$ where $f_{i,\alpha}\in\ZZ[\und{a}]$ for any $i=1,...,s$ and $\alpha\in\NN^{m}$. Assume that $K$ is an algebraically closed field and $\und{a}\in K^{n}$. Then $f_{i,\alpha}(\und{a})\in K$ for any $i=1,...,s$, $\alpha\in\NN^{m}$ and thus it follows from Hilbert's Nullstellensatz that $\varphi(\und{a})$ holds if and only if $1\notin\langle F_{1},...,F_{s}\rangle$. Corollary 2.3 implies that $1\notin\langle F_{1},...,F_{s}\rangle$ is equivalent with non-existence of polynomials $H_{1},...,H_{s}\in K[x_{1},...,x_{m}]$ such that $1=H_{1}F_{1}+...+H_{s}F_{s}$ and $\deg(H_{i})\leq\gamma_{m,d}(0)=d'$ for any $i=1,...,s$. The fact that $\deg(H_{i})\leq d'$ enables to write the latter condition in the first order language.

We introduce some sets of variables. Assume that $C_{i}=\{c_{i,\beta}\}_{|\beta|\leq d'}$ where $\beta\in\NN^{m}$ and $i=1,...,s$. Let $H_{i}\in\ZZ[C_{i}][\und{x}]$ be a polynomial of the form $H_{i}=\sum_{|\beta|\leq d'}c_{i,\beta}\cdot\und{x}^{\beta}$ for $i=1,...,s$. Set $C=\bigcup_{i=1}^{s}C_{i}$ and consider the formula $\psi(\und{a})=\forall_{C}H_{1}F_{1}+...+H_{s}F_{s}\neq 1$ which is equivalent with $\varphi(\und{a})$. Observe that $$H_{1}F_{1}+...+H_{s}F_{s}=\sum_{|\delta|\leq d+d'}(\sum_{\beta+\alpha=\delta}c_{1,\beta}f_{1,\alpha}+...+c_{s,\beta}f_{s,\alpha})\und{x}^{\delta}$$ where $\delta\in\NN^{m}$, and hence the formula $\psi(\und{a})$ expresses the non-existence of solution of some system of linear equations with the set $C$ as a set of variables. This system can be written in such a way that the matrices $A=A^{d,d'}_{F_{1},...,F_{s}}$ and $\widehat{A}=\widehat{A}^{d,d'}_{F_{1},...,F_{s}}$ are its coefficient matrix and augmented matrix, respectively. Then it follows from the Kronecker-Capelli theorem that $\psi(\und{a})$ holds if and only if $\rk(\widehat{A})>\rk(A)$ where $\rk(M)$ denotes the rank of the matrix $M$. This is equivalent with $\Delta^{d,d'}_{F_{1},...,F_{s}}(\und{a})$. \epv

As a direct consequence of our considerations we get the following computable condition for the existence of a common root of multivariate polynomials.  

\begin{cor} Assume that $K$ is an algebraically closed field, $d$ is a natural number, $F_{1},...,F_{s}\in K[x_{1},...,x_{m}]$ and $F_{i}=\sum_{|\alpha|\leq d}a_{i,\alpha}\cdot\und{x}^{\alpha}$ for $i=1,...,s$. Set $d'=\gamma_{m,d}(0)$. The polynomials $F_{1},...,F_{s}$ have a common root if and only if $\rk(\widehat{A})>\rk(A)$ where $A$ and $\widehat{A}$ are matrices obtained from $A^{d,d'}_{F_{1},...,F_{s}}$ and $\widehat{A}^{d,d'}_{F_{1},...,F_{s}}$, respectively, by replacing the elements $f_{i,\alpha}$ by $a_{i,\alpha}$ for any $i=1,...,s$, $\alpha\in\NN^{m}$. \epv
\end{cor}

{\bf Proof.} The proof is a simplified version of the proof of Theorem 3.2. \epv

\section{Applications}

In this section we present some applications of Theorem 3.2 in mathematics and physics, especially in quantum information theory. We concentrate on the problem of the existence of common invariant subspaces of square complex matrices and related problems. In that sense, we continue our research (and generalize the results) from \cite{JP} and \cite{PJ}, see also \cite{JKP} and \cite{PKJ}.

We give this section an expository character and leave the details for further papers. We recommend \cite{BZ} (see also \cite{HZ}) as a comprehensive monograph on quantum information theory and quantum mechanics in general.

Assume that $A,A_{1},...,A_{s}$ are $n\times n$ matrices over the field $\mathbb{C}$ of complex numbers and $V$ is a subspace of $\mathbb{C}^{n}$. We say that $V$ is \textit{$A$-invariant} if and only if $Av\in V$ for any $v\in V$. We say that $V$ is a \textit{common invariant subspace} of $A_{1},...,A_{s}$ if and only if $V$ is $A_{i}$-invariant for any $i=1,...,s$.

The problem of the existence of common invariant subspaces of square complex matrices appears in many areas of mathematics and physics. Therefore, computable conditions for the existence of such subspaces are heavily studied. In \cite{Sh} the author gives a computable condition for the existence of a common eigenvector (i.e. a common invariant subspace of dimension one) of two matrices. This result is generalized to a finite number of matrices in \cite{JP}, see also \cite{PJ} for similar concepts. In \cite{AGI}, \cite{AI}, \cite{GI} and \cite{Ts} only two matrices are considered, but the authors study common invariant subspaces of dimensions higher than one. In this case it is often assumed that given matrices have pairwise different eigenvalues. This assumption is made in \cite{GI} and \cite{Ts} where the authors reduce the general problem to the question of the existence of a common eigenvector of suitable \textit{compound matrices}, see \cite{MM}.

The general version of the problem, with arbitrary finite number of matrices and arbitrary dimension of common invariant subspaces, was solved only in 2004 in \cite{ArPe}. In the solution some basic techniques of Gr{\"o}bner bases theory and algebraic geometry are used.

Here we apply Theorem 3.2 (or Corollary 3.3) to solve the general problem of the existence of a common invariant subspace. Assume that $A_{1},...,A_{s}$ are complex $n\times n$ matrices. Let $V=\{v_{i}^{j}|i=1,...,k,j=1,...,n\}$ be a set of variables and $\widehat{V}$ the set of all $\mathbb{C}$-linear combinations of elements of $V$. Set $\und{v}_{i}=\left[\begin{matrix}v_{i}^{1}&&\hdots&&v_{i}^{n}\end{matrix}\right]^{T}$ for $i=1,...,k$ and denote by $M_{V}$ the augmented matrix $(\und{v}_{1}|...|\und{v}_{k})$. The formula $\rk(M_{V})=k$ states that the vectors $\und{v}_{1},...,\und{v}_{k}$ are linearly independent, and can be written in the first order language. The formula $\bigwedge_{j=1}^{k}A_{i}\und{v}_{j}\in\widehat{V}$, for any $i=1,...,s$, states that $\widehat{V}$ is $A_{i}$-invariant, and can be written in the first order language. Thus the first order formula $$\exists_{V}\rk(M_{V})\wedge(\bigwedge_{i=1}^{s}\bigwedge_{j=1}^{k}A_{i}\und{v}_{j}\in\widehat{V}),$$ where $k\leq n$, expresses the existence of a common invariant subspace of $A_{1},...,A_{s}$ of dimension $k$. It is easy to see that this is in fact a common root formula and hence Theorem 3.2 (or Corollary 3.3) yields its equivalent quantifier-free form. We call this quantifier-free form a \textit{$\CIS_{k}$-formula} for $A_{1},...,A_{s}$. Such a formula can be viewed as a computable condition for the existence of a common invariant subspace of $A_{1},...,A_{s}$ of dimension $k$.

Common invariant subspaces, sometimes satisfying some additional conditions, play a prominent role in quantum information theory. We show this role on two examples concerning quantum channels (so our treatment of the subject is far from being complete): irreducible quantum channels and decoherence-free subspaces. In these examples we apply $\CIS_{k}$-formulas and Theorem 3.2 to generalize some results from \cite{JP} and \cite{PJ}.

Assume that $\mathbb{M}_{n}(\mathbb{C})$ is the vector space of all $n\times n$ complex matrices. A \textit{quantum channel} is a trace preserving completely positive map $\Phi:\mathbb{M}_{n}(\mathbb{C})\ra\mathbb{M}_{n}(\mathbb{C})$ (we refer to \cite{HZ} for all the definitions). It follows from \cite[5.2.3]{HZ} that there are matrices $A_{1},...,A_{s}\in\mathbb{M}_{n}(\mathbb{C})$ such that $\Phi(X)=\sum_{i=1}^{s}A_{i}XA_{i}^{*}$ for any $X\in\mathbb{M}_{n}(\mathbb{C})$ where $A^{*}$ denotes the matrix adjoint to $A$. 

Important subclass of the class of all quantum channels is formed by \textit{irreducible} quantum channels. We refer to \cite{BZ} and \cite{HZ} for the definition and main properties of these channels. It is proved in \cite{Fa} that a quantum channel $\Phi(X)=\sum_{i=1}^{s}A_{i}XA_{i}^{*}$ is irreducible if and only if the matrices $A_{1},...,A_{s}$ do not have a nontrivial common invariant subspace. Hence the $\CIS_{k}$-formulas for $A_{1},...,A_{s}$ provide a computable condition for irreducibility of $\Phi$. This generalizes the main results of \cite{JP}, see especially Sections 3 and 4 of \cite{JP}.

Quantum channels are used to transmit quantum information. Unfortunately, quantum information may be easily corrupted by a number of factors, see \cite{BL}. Any such a factor is described as a \textit{decoherence}. A way to overcome the effects of decoherence is to "hide" quantum information from the environment in some "quiet corner". This quiet corner is called the \textit{decoherence-free subspace} (DFS). 

There are few different mathematical definitions of DFS in the literature, see \cite{KMSW} for the details. In \cite{PJ} we define DFS as the \textit{common reducing unitary subspace}. We recall this definition below.

Assume that $A,A_{1},...,A_{s}\in\mathbb{M}_{n}(\mathbb{C})$ and $W$ is a subspace of $\mathbb{C}$.  We say that $W$ is a \textit{reducing} subspace of $A$ (or \textit{$A$-reducing}) if and only if $W$ is an invariant subspace for $A$ and $A^{*}$. We say that $W$ is a \textit{common reducing subspace} of $A_{1},...,A_{s}$ if and only if $W$ is $A_{i}$-reducing for any $i=1,...,s$. 

Assume that $A_{1},...,A_{s}\in\mathbb{M}_{n}(\mathbb{C})$ and $\Phi(X)=\sum_{i=1}^{s}A_{i}XA_{i}^{*}$ is a quantum channel. A nonzero subspace $W$ of $\mathbb{C}^{n}$ is a \textit{common reducing unitary subspace} (or a \textit{decoherence-free subspace}) for $\Phi$ if and only if $W$ is a common reducing subspace of $A_{1},...,A_{s}$ and there exists a unitary matrix $U\in\mathbb{M}_{n}(\mathbb{C})$ and complex numbers $g_{1},...,g_{s}$ such that $A_{i}w=(g_{i}U)w$ for any $w\in W$ and $i=1,...,s$. 

The conditions that $U\in\mathbb{M}_{n}(\mathbb{C})$ is a unitary matrix and $A_{i}w=(g_{i}U)w$ for any $w\in W$ and $i=1,...,s$ can be written in the first order language. Hence there is a formula expressing the existence of a common reducing unitary subspace of dimension $k$. This formula is similar to $\CIS_{k}$-formula. Consequently, Theorem 3.2 provides a computable condition for the existence of decoherence-free subspaces. This generalizes the main results of \cite{PJ}, see especially Section 3 of \cite{PJ}.

The contents of the section reveal that there is an impact of quantifier elimination theory on applied mathematics. This impact has been recently noticed in a number of papers, see for example \cite{WCP-G} and \cite{SZ}.

The results of Section 3 imply that every problem which can be written in the first order language of fields can be equivalently expressed as a computable condition. Moreover, Theorem 3.2 provides the exact form of this condition. It is our opinion that this opens the possibility for other applications of quantifier elimination in mathematical sciences.

\section*{Acknowledgements} This research has been supported by grant No. DEC-2011/02/A/ST1/00208 of National Science Center of Poland.

\noindent Grzegorz Pastuszak\\Faculty of Mathematics and Computer Science\\ Nicholaus
Copernicus University\\ Chopina 12/18\\ 87-100 Toru\'n, Poland\\
past@mat.uni.torun.pl

\end{document}